\documentclass{article}
\usepackage{amsmath}
\usepackage{graphicx}
\usepackage{latexsym}
\usepackage{hyperref}
\usepackage{fancyhdr}
\usepackage{amsfonts}
\usepackage[all]{xy}
\newcommand{\R}{\mathbb{R}}

\def\p{\partial}

\def\b{\beta}

\newcommand{\ri}{\rightarrow}
\newcommand{\na}{\nabla}
\newcommand{\pa}{\partial}
\newcommand{\no}{\nonumber}

\newcommand{\ra}{\rangle}
\newcommand{\fr}{\frac}

\newcommand\be{\begin{eqnarray}}
\newcommand\en{\end{eqnarray}}
\def\en{\end{eqnarray}}

\def\inpo{\int_{\partial\Omega}}
\def\D{\Delta}

\def\la{\lambda}

\def\om{\Omega}

\def\({\left(}
\newcommand{\Om}{\Omega}
\newcommand{\laa}{\lambda}
\newcommand{\ov}{\overline}

\def\){\right)}

\newtheorem{theorem}{Theorem}[section]
\newtheorem{lemma}[theorem]{Lemma}

\newtheorem{conjecture}[theorem]{Conjecture}

\begin{document}
\setcounter{page}{1}
\title{Isoperimetric Bounds  for Lower Order Eigenvalues}
\renewcommand{\thefootnote}{}
\footnotetext{

\hspace*{2ex}Key words and phrases: Isoperimetric inequalities, eigenvalues, Laplacian, biharmonic Steklov
problems, Wentzell-Laplace operator.}
\author{Fuquan Fang and Changyu Xia}
\date{}
\maketitle ~~~\\[-15mm]

\begin{abstract}
New isoperimetric inequalities for lower order
eigenvalues of the Laplacian on closed
hypersurfaces,  of the biharmonic Steklov problems and of the Wentzell-Laplace on bounded domains in a Euclidean space are proven. Some open questions for further study are also proposed.
 \end{abstract}

\markright{\hfill Fuquan Fang and Changyu Xia \hfill}

\section{Introduction and the main results}
\renewcommand{\thesection}{\arabic{section}}
\renewcommand{\theequation}{\thesection.\arabic{equation}}
\setcounter{equation}{0} \label{intro}

Let $(M, g)$ be a  closed Riemannian manifold of dimension $\geq 2$. The spectrum of the Laplace operator on $M$ provides a sequence of global Riemannian invariants
\be\no
0=\la_0(M)<\la_1(M)\leq \la_2(M)\leq\cdots\nearrow\infty.
\en
We adopt the convention that each eigenvalue is
repeated according to its multiplicity. An important issue in spectral geometry is to obtain good estimates for  these and other eigenvalues in terms of  the   geometric data of the manifold $M$ such as the volume, the diameter, the
curvature, the  isoperimetric constants, etc.   See \cite{be},\cite{ber},\cite{c2},\cite{ds},\cite{sy} for  references.

On the other hand, after the seminal works of Bleecker-Weiner \cite{bw} and Reilly \cite{r1}, the following approach is developed:  the manifold $(M, g)$ is immersed
isometrically into another Riemannian manifold. One then gets good estimates for $\lambda_k(M)$, mostly for $\lambda_1(M)$, in termos of the extrinsic  geometric quantities of $M$.  See for example \cite{bw}, \cite{dmwx}, \cite{ei}, \cite{gr}, \cite{h}, \cite{wx1}, \cite{x1}.
Especially relevant for us is the quoted work of
 Reilly \cite{r1}, where  he  obtained the following remarkable isoperimetric inequality for the first positive eigenvalue $\la_1(M)$ in the case that $M$ is embedded as a
hypersurface bounding a domain $\om$ in $\mathbb{R}^{n}$:
\be\label{0.1}
\la_1(M)\leq\fr{n-1}{n^2}\cdot\fr{|M|^2}{|\om|^2}.
\en
Here $|M|$ and $|\om|$ denote the Riemannian $(n-1)$-volume of $M$ and the Riemannian $n$-volume of $\om$, respectively. Moreover, equality holds in (\ref{0.1}) if and only if $M$ is  a
round sphere.   Our first result improves  (\ref{0.1}) to the sum of the first $n$ non-zero eigenvalues of the Laplace operator on
$M$.
\begin{theorem} \label{th1} Let $M$ be a closed embedded hypersurface  bounding a domain $\om$ in $\mathbb{R}^{n}$. Then the first $n$ non-zero eigenvalues of the Laplacian on
$M$ satisfy
\be\label{0.2}
\sum_{i=1}^{n}\la_i\leq \fr{n-1}{n}\cdot\fr{|M|^2}{|\om|^2}
\en
and
\be\label{0.3}
\sum_{i=1}^{n}\la_i\leq \fr{(n-1)\sqrt{|M|}}{|\om|}\left(\int_M H^2\right)^{1/2},
\en
where $H$ stands for the mean curvature of $M$. Moreover, equality holds in either of (\ref{0.2}) and (\ref{0.3}) if and only if $M$ is  a sphere.
\end{theorem}

In the second part of this paper we study eigenvalues of fourth order Steklov problems. Let
$\om$ be an $n$-dimensional compact Riemannian manifold with boundary and  $\Delta$ and $\ov{\Delta}$
be the Laplace
operators on $\om$ and $\pa\om$, respectively. Consider
the eigenvalue problem
\be\label{0.4}
 \left\{\begin{array}{l}
 \Delta^2 u= 0  \ \ {\rm in \ \ } \om, \\
\pa_{\nu}u = \pa_{\nu}(\D u)+\xi u =0 \ \ {\rm on \ \ } \pa\om,
\end{array}\right.
\en
where $\pa_{\nu}$ denotes the outward unit normal derivative. This problem was first
discussed by J. R. Kuttler and V. G. Sigillito  \cite{ks} in the case where $\om$ is a bounded domain in $\R^n$.   The eigenvalue problem (\ref{0.4})
is important in biharmonic analysis and elastic mechanics. In the two dimensional case, it describes the deformation $u$ of the linear elastic
supported plate $\Om$ under the action of the transversal exterior
force $f(x) = 0, \ x\in\Om $ with Neumann boundary condition
$\pa_{\nu}u|_{\pa \Om}=0$ (see, \cite{tg},\cite{v},\cite{xw2}).  In addition, the first nonzero eigenvalue $\xi_1$ arises as an optimal constant in an a priori inequality (see \cite{ks}).
The eigenvalues of
the problem (\ref{0.4}) form a discrete and increasing sequence (counted with
multiplicity):
\be
0=\xi_0<\xi_1\leq\xi_2\leq\cdots\nearrow +\infty.
\en
Let $\mathcal{D}_k$ be the space of harmonic homogeneous polynomials
in $\R^n$ of degree $k$ and denote by $\mu_k$ the dimension
of $\mathcal{D}_k$, $k=0, 1,\cdots,$. For the $n$-dimensional Euclidean ball  with radius $R$, the eigenvalues of (\ref{0.4}) are $\xi_k=k^2(n+2k)/R^3, k=0, 1, 2, \cdots, $ and the multiplicity of $\xi_k$ is $\mu_k$ (see \cite{xw2}, Theorem 1.5 ). When $\om$ has nonnegative Ricci curvature with strictly convex boundary, a lower bound for
$\xi_1(\om)$ has been given in \cite{xw1}.  On the other hand,   an isoperimetric upper bound  for $\xi_1(\om)$ has been proven  for the case where $\om$ is a bounded domain in $\R^n$ (see \cite{xw2}, Theorem 1.6 ). We have an isoperimetric inequality for the sum of the reciprocals of the first $n$ nonzero eigenvalues of
the problem (\ref{0.4}) on  bounded domains in $\R^n$.
\begin{theorem} \label{th3} Let $\Om$ be a bounded domain with smooth boundary in $\mathbb{R}^n$. Then the first $n$ nonzero eigenvalues of
the problem (\ref{0.4}) satisfy
\be\label{0.5}
\sum_{i=1}^n\fr 1{\xi_i}\geq \fr{n^2|\om|\left(\fr{|\om|}{\omega_n}\right)^{2/n}}{(n+2)|\pa\om|},
\en
with equality holding if and only if $\om$ is a ball, where $\omega_n$ denotes the volume of the unit ball in
$\R^n$.
\end{theorem}
Now we come to another Steklov problem for the bi-harmonic operator. Let $\om$  be a bounded domain in $\R^{n}$ and $\tau$ a
positive constant. Denote by $\nabla^2$ and  $\nabla$ the Hessian on $\R^n$ and the gradient operator on $\om$, respectively. Consider the following Steklov problem of fourth
 order
\begin{eqnarray}\label{0.6}
\left\{\begin{array}{ccc} \D^2 u-\tau \D u
=0,&&~\mbox{in} ~~ \om, \\[2mm]
\frac{\p^2u}{\p\nu^2}=0, &&~~\mbox{on}~~\partial \om,\\[2mm]
\tau\frac{\p u}{\p \nu}-\mathrm{div}_{\p M}\(\nabla^2
u(\nu)\)-\frac{\p\D u}{\p\nu}=\lambda u,
&&~~\mbox{on}~~\partial \om.
\end{array}\right.
\end{eqnarray}
This problem has a discrete spectrum which can be listed as
$$0=\lambda_{0,\tau}<\lambda_{1,\tau}\leq\lambda_{2,\tau}\leq\cdots\leq\lambda_{k,\tau}\nearrow +\infty.$$
The eigenvalue $0$ is simple and the corresponding eigenfunctions
are constants. Let $u_0, u_1,..., u_{k},\cdots, $ be  the  eigenfunctions of problem (\ref{0.6}) corresponding to the eigenvalues  $0=\lambda_{0,\tau},\ \lambda_{1,\tau},\cdots, \lambda_{k,\tau}, \cdots, $.  For each $k=1,\cdots,$ we have the following variational characterization
 \begin{small}
 \begin{eqnarray}\label{0.7}
\lambda_{k,\tau}=\mathrm{min}\left\{\frac{\int_{\om}\(|\na^2u|^2+\tau|\na
u|^2\)}{\int_{\pa\om} u^2}\Bigg| u\in H^2(\om), u\neq0,
\int_{\pa\om} u u_j=0, j=0,\cdots,k-1\right\}.
\end{eqnarray}
\end{small}

The eigenvalues and eigenfunctions on the ball   in $\R^n$ have been determined by  Buoso-Provenzano in \cite{bp}. In particular, if $\mathbf{B}_R^n$ is the ball of radius $R$ centered at the origin in $\R^n$,
then  \be
\lambda_{1,\tau}(\mathbf{B}_R^n)=\lambda_{2,\tau}(\mathbf{B}_R^n)=\cdots=\lambda_{n,\tau}(\mathbf{B}_R^n)=\fr{\tau}{R}
\en
 and the corresponding eigenspace is generated by $\{x_1,\cdots, x_n\}$.  Buoso and Provenzano \cite{bp} also proved the following
 isoperimetric inequality for the sums of the
reciprocals of the first $n$ non-zero eigenvalues:
\be\label{bp}
\sum_{i=1}^n\fr 1{\lambda_{1,\tau}(\om)}\geq \fr n{\tau}\left(\fr{|\om|}{\omega_n}\right)^{1/n}
\en
with equality holding if and only if $\om$ is a ball.
  Further study for the eigenvalues of the problem (\ref{0.6}) has been made in \cite{bcp},
\cite{dmwx}, \cite{xw2}, etc.  Our next result is an isoperimetric inequality for the sum of
the first $n$ non-zero eigenvalues of the problem (\ref{0.6}).

\begin{theorem}\label{th4} Let $\Omega$ be a bounded domain with smooth boundary $\p \Omega$ in $\R^n$. Denoting by
$\lambda_{i,\tau}$ the $i$-th eigenvalue of the (\ref{0.6}),  we
have
\begin{eqnarray}\label{th4.1}
\sum_{j=1}^{n}\lambda_{j,\tau}\leq
\frac{\tau|\pa\Omega|}{|\om|}.
\en
Equality holds in (\ref{th4.1}) if and only if $\Omega$ is a ball.
\end{theorem}

The final part of the present paper concerns   the eigenvalue
problem with Wentzell boundary conditions:
\be\label{int1}\left\{\begin{array}{l} \Delta u =0 \ \ \ \ \ \
\ \ \ \ \ \ \ \ \ \ \ \ \ \ \ \ {\rm in  \ }\ \om,\\ -\beta\ov{\Delta}
u+\pa_{\nu} u= \laa u\ \ \ \quad~ {\rm on  \ } \pa \om,
\end{array}\right.
\en where $\beta$ is a nonnegative constant, $\om$  is a compact Riemannian manifold of dimension $n\geq 2$ with non-empty boundary,   $\Delta$  and $\ov{\Delta}$ denote  the the Laplacian on $\om$ and $\pa\om$, respectively. When $\beta=0$, (\ref{int1}) becomes the  Steklov problem:
\be\label{int2}\left\{\begin{array}{l} \Delta u =0 \ \ \ \ \ \
\ \ \   {\rm in  \ }\ \om,\\ \pa_{\nu} u= p u\ \ \ \quad~ {\rm on  \ } \pa \om,
\end{array}\right.
\en
which has been studied  extensively ( see \cite{b}, \cite{bpr},\cite{ceg},\cite{e1}-\cite{fs3},\cite{hps},\cite{ks},\cite{st}, \cite{wx2},\cite{xw2},\cite{xi} ). The spectrum of the problem (\ref{int1})  consists in an increasing sequence
 \be \no\laa_{0, \beta}=0<\laa_{1, \beta}\leq\laa_{2,
\beta}\leq\cdots\nearrow +\infty, \en with corresponding real
orthonormal (in $L^2(\pa \om)$ sense) eigenfunctions $u_0, u_1,
u_2,\cdots.$
Consider the Hilbert space
\be H(\om)=\{ u\in H^1(\om), {\rm  Tr}_{\pa \om}(u)\in H^1(\pa \om) \}, \en
where ${\rm Tr}_{\pa \om}$ is the trace operator. We define on $H(\om)$
the two bilinear forms \be A_{\beta}(u,
v)=\int_{\om}\nabla u\cdot\nabla v +\beta\int_{\pa \om} \ov{\na}
u\cdot\ov{\na} v, \ B(u, v)=\int_{\pa \om}uv, \en where,
$\nabla $ and $\ov{\na}$ are the gradient operators on $\om$ and $\pa
\om$, respectively. Since  we assume  that $\beta$  is nonnegative,
the two bilinear forms are positive and the variational
characterization for the $k$-th eigenvalue is \be \label{vc}\laa_{k,
\beta}=\min\left\{\frac{A_{\beta}(u, u)}{B(u, u)}, u\in H(\om), u\neq
0,\ \int_{\pa \om} u u_i=0, i=0,\cdots,k-1\right\}. \en When $k=1$,
the minimum is taken over the functions orthogonal to the
eigenfunctions associated to $\laa_{0, \beta}= 0,$ i.e., constant
functions.

If $\om=\mathbf{B}_R^n$,  then
\cite{dkl}
$$\laa_{1, \beta}=\laa_{2, \beta}=\cdots \laa_{n, \beta}=\fr{(n-1)\beta + R}{R^2}$$
and  the corresponding eigenspace is generated by $\{x_i,  i=1,\cdots,n\}$. For the Steklov problem (\ref{int2}), Brock \cite{b} showed that if $\om$ is a bounded domain with smooth boundary in $\R^n$,
then  the first $n$ nonzero eigenvalues of $\om$ satisfy
\be\label{th0.0}
\sum_{i=1}^n\fr 1{p_i(\om)} \geq
n\left(\fr{|\om|}{\omega_{n}}\right)^{\fr 1n}, \en
with equality holding if and only if $\om$ is a ball.
Brock's theorem has been generalized to the  eigenvalues of the
problem  (\ref{int1}) in \cite{dmwx}. We prove

\begin{theorem} \label{th7} Let $\beta\geq 0$ and $\Omega$ be a bounded domain with smooth boundary $\p \Omega$ in
 $\R^n$.  Denote by $\laa_{1, \beta}\leq\laa_{2, \beta}\leq\cdots\leq \laa_{n, \beta}$
 the first $n$ non-zero
 eigenvalues of the following problem with the Wentzell boundary condition.
\be\label{th7.1}\left\{\begin{array}{l} \Delta u =0 \ \ \ \ \ \
\ \ \ \ \ \ \ \ \ \ \ \ \ \ \ \ \ {\rm in  \ }\ \om,\\ -\beta\ov{\Delta}
u+\pa_{\nu} u= \laa u\ \ \ \ \quad ~~ {\rm on  \ } \pa \om,
\end{array}\right.
\en
Then we have
\be\label{th7.2}
\sum_{i=1}^n \lambda_{i,\b}\leq \fr{|\pa\om|}{|\om|}+\fr{(n-1)\beta}n\cdot\fr{|\pa\om|^2}{|\om|^2}.
 \en
Furthermore, equality holds in (\ref{th7.2}) if and only if $\om$
is a ball.
 \end{theorem}
Taking $\beta=0$ in (\ref{th7.2}), we have a new isoperimetric inequality for the first $n$ nonzero Steklov eigenvalues of a bounded domain $\om\subset\R^n$:
\be\label{th5.3}
\sum_{i=1}^n p_i(\om)\leq \fr{|\pa\om|}{|\om|}
\en
with quality holding  if and only if $\om$
is a ball.

It has been conjectured by Henrot \cite{he1} that the first $n$ nonzero Steklov eigenvalues of a bounded domain $\om\subset\R^n$ satisfy
\be\label{th5.3}
\overset{n}{\underset{i=1}{\prod}} p_i(\om)\leq\fr{\omega_n}{|\om|}
\en
which is stronger than Brock's inequality (\ref{th0.0}). If $n=2$, or $n\geq 3$ and $\om$ is convex, then (\ref{th5.3}) is true ( see \cite{hps}, \cite{he2} ). This result can be extended to eigenvalues of the problem (\ref{int1}). Namely, we have
\begin{theorem} \label{th8} Let the notation be as in Theorem \ref{th7} and when $n\geq 3$,  assume
further that $\om$ is convex. Then
\be\label{th8.1}
\overset{n}{\underset{i=1}{\prod}} \la_{i,\beta}\leq\left(1+\fr{(n-1)\beta|\pa\om|}{n|\om|}\right)^n\cdot\fr{\omega_n}{|\om|}
\en
 with quality holding  if and only if $\om$
is a ball.
\end{theorem}

\section{A Proof of Theorem \ref{th1}}
\setcounter{equation}{0} In this section, we give a
\vskip0.3cm
{\it Proof of Theorem \ref{th1}}.
  Let $\D$ and  $\ov{\D}$ be the Laplace operators on $\R^n$ and $M$, respectively,  and let $\{u_i\}_{i=0}^{+\infty}$
be an orthonormal system of eigenfunctions  corresponding to the
eigenvalues \be\label{pth3.1} 0=\la_0<\la_1\leq
\la_2\leq\cdots\ri\infty \en of the Laplacian of $M$, that is,
\be\label{th3.2} \ov{\Delta} u_i = -\la_i u_i, \quad  \int_M u_iu_j
=\delta_{ij}. \en We have $u_0=1/\sqrt{|M|}$ and for each
$i=1,\cdots,$ the Rayleigh-Ritz characterization for $\la_i$ is given by
\be \label{pth1.1} \la_i=\underset{u\neq 0, \int_M
uu_j=0, j=0,\cdots i-1}{\min}\fr{\int_M |\ov{\na} u|^2}{\int_M u^2},\en
being $\ov{\na}$ the gradient operator on $M$.

In order to obtain good upper bound for $\la_i$, we need to choose nice trial functions
$\phi_i$ for each of the eigenfunctions $u_i$ and insure that these
are orthogonal to the preceding eigenfunctions $u_0,\cdots,
u_{i-1}$. We note that the coordinate functions are eigenfunctions corresponding to the
first eigenvalue of the hypersphere in $\R^n$. For the $n$ trial functions $\phi_1, \phi_2, \cdots,
\phi_n,$ we simply choose the $n$ coordinate functions: \be
\phi_i=x_i, \ \ {\rm for}\ \ i=1,\cdots, n, \en but before we can
use these we need to make adjustments  so that $\phi_i \perp{\rm
span}\{u_0,\cdots, u_{i-1}\}$ in $L^2(\pa\om)$. By
translating the origin appropriately we can assume that \be
\int_{M} x_i=0, \ i=1,\cdots,n, \en that is, $x_i\perp u_0$.

Nextly we show that a rotation of the axes can be made  so that  \be\label{pth1.} \int_{M}\phi_j u_i=\int_{M}
x_j u_i=0, \en for $j=2,3,\cdots, n$ and $i=1,\cdots,j-1$. In fact,
 let us  define an $n \times n$ matrix $P=\(p_{ji}\),$ where
$p_{ji}=\int_{M} x_j u_i$, for $i,j=1,2,\cdots,n.$ Using the
orthogonalization of Gram and Schmidt (QR-factorization theorem), one
can find  an upper triangle matrix $T=(T_{ji})$ and an
orthogonal matrix $U=(a_{ji})$ such that $T=UP$, that is,
\begin{eqnarray*}
T_{ji}=\sum_{k=1}^n a_{jk}p_{ki}=\int_M \sum_{k=1}^n  a_{jk}x_k  u_i
=0,\ \  1\leq i<j\leq n.
\end{eqnarray*}
Letting $y_j=\sum_{k=1}^n  a_{jk}x_k$, we have
\begin{eqnarray}\label{c2}
\int_M y_j  u_i =\int_M \sum_{k=1}^n a_{jk}x_k u_i
=0,\ \  1\leq i<j\leq n.
\end{eqnarray}
Since $U$ is an orthogonal matrix,  $y_1, y_2, \cdots, y_n$ are
also coordinate functions on $\R^n$. Thus, denoting these
coordinate functions still by $x_1, x_2,\cdots, x_n$, one arrives at
the  condition (\ref{pth1.}).
It  follows from (\ref{pth1.1}) that
\be\label{pt1.4} \la_i\int_M  x_i^2\leq \int_M |\ov{\na}
x_i|^2, \ i=1,\cdots, n,
\end{eqnarray}
with equality holding if and only if \be\label{pth1.5} \ov{\D} x_i=
-\la_i x_i. \en
  Integrating the equality
\be\label{pth1.2}
\fr 12\Delta x_i^2=1
\en
on $\om$ and using the divergence theorem, one gets
\be\label{pth1.3}
|\om|=\int_M x_i\pa_{\nu}x_i, \  i=1,\cdots, n,
\en
where $\nu$ denotes the outward unit normal of $\pa\om=M$. Taking the square of (\ref{pth1.3}) and using the H\"older inequality, we infer
\be\label{pth1.6}
|\om|^2\leq\left(\int_M x_i^2\right)\left(\int_M(\pa_{\nu}x_i)^2\right), \  i=1,\cdots, n.
\en
Multiplying (\ref{pt1.4}) by $\int_M(\pa_{\nu}x_i)^2$ and using (\ref{pth1.6}), we have
\be\label{pth1.7}
|\om|^2\la_i\leq \left(\int_M |\ov{\na}x_i|^2\right)\left(\int_M(\pa_{\nu}x_i)^2\right) , \ i=1,\cdots, n.
\en
Observing  that on $M$
\be
1=|\na x_i|^2=|\ov{\na}x_i|^2+(\pa_{\nu}x_i)^2,
\en
one deduces from (\ref{pth1.7}) that
\be\label{pth1.8}
|\om|^2\la_i\leq \left(|M|-\int_M (\pa_{\nu}x_i)^2\right)\left(\int_M(\pa_{\nu}x_i)^2\right), \ i=1,\cdots, n.
\en
Summing over $i$ and using Cauchy-Schwarz inequality, we get
\be\label{pth1.9}
|\om|^2\sum_{i=1}^{n}\la_i&\leq& |M|\sum_{i=1}^n\left(\int_M(\pa_{\nu}x_i)^2\right)-\sum_{i=1}^n\left(\int_M (\pa_{\nu}x_i)^2\right)^2\\ \no &\leq& |M|^2-\fr 1n\left(\sum_{i=1}^n\int_M (\pa_{\nu}x_i)^2\right)^2\\ \no
&=& \fr{(n-1)|M|^2}n.
\en
This proves (\ref{0.2}).

To prove (\ref{0.3}), we use  divergence theorem and  H\"older inequality to get
\be\label{pth1.10}
\int_M |\ov{\na}x_i|^2=-\int_M x_i\ov{\D}x_i\leq \left(\int_M x_i^2\right)^{1/2}\left(\int_M (\ov{\D}x_i)^2\right)^{1/2},
\en
which, combining with (\ref{pt1.4}), gives
\be\label{pth1.11}
\la_i\left(\int_M x_i^2\right)^{1/2}\leq\left(\int_M (\ov{\D}x_i)^2\right)^{1/2}.
\en
Multiplying (\ref{pth1.11}) by $\left(\int_M (\pa_{\nu}x_i)^2\right)^{1/2}$ and using (\ref{pth1.6}), we have
\be\label{pth1.12}
\la_i|\om|\leq \left(\int_M (\ov{\D}x_i)^2\right)^{1/2}\left(\int_M (\pa_{\nu}x_i)^2\right)^{1/2}, \ i=1,\cdots,n.
\en
Summing over $i$ and using Cauchy-Schwarz inequaty,
we infer
\be\label{pth1.13}
|\om|\sum_{i=1}^n\la_i&\leq& \left(\sum_{i=1}^n\int_M (\ov{\D}x_i)^2\right)^{1/2}\left(\sum_{i=1}^n\int_M (\pa_{\nu}x_i)^2\right)^{1/2}\\ \no &=&
(n-1)\left(\int_M H^2\right)^{1/2}|M|^{1/2},
\en
where, in the last equality, we have used the fact that
\be\label{pth1.14} \ov{\Delta} x \equiv (\ov{\Delta} x_1,\cdots, \ov{\Delta} x_n)
=(n-1){\mathbf H}, \en being ${\mathbf H}$  the mean curvature
vector of $M$ in ${\R}^n$. Hence, (\ref{0.3}) holds.

If the equality holds in (\ref{0.2}), then the inequalities
(\ref{pt1.4}), (\ref{pth1.6}), (\ref{pth1.8}) and (\ref{pth1.9}) must take equality sign. It then follows that (\ref{pth1.5}) holds,
\be \int_M (\pa_{\nu}x_1)^2=\int_M (\pa_{\nu}x_2)^2=\cdots=\int_M (\pa_{\nu}x_n)^2,
\en
and so
\be
\la_1=\la_2=\cdots=\la_n.
\en
Thus, the
position vector $x=(x_1,\cdots, x_n)$ when restricted on $\pa\Om$
satisfies
\be\label{p1} \Delta x = -\la_1
(x_1,\cdots, x_n).\en  Combining (\ref{p1}) and
(\ref{pth1.14}), we have \be\label{p3} x=-\fr{(n-1)}{\la_1}{\mathbf H}, \ \
\ \ {\rm on\ \ \ }M. \en Consider the function $h = |x|^2 :
M\rightarrow {\R}$. It is easy to see from (\ref{p3}) that \be\no Z
h =2\langle Z, x\rangle = 0, \ \ \ \forall Z\in
{\mathfrak{X}}(M). \en Thus $g$ is a constant function and so
$M$ is a hypersphere. If
the equality holds in (\ref{0.3}), one can use similar arguments to deduce that $M$ is a hypersphere in $\R^n$.  \hfill
$\Box$
\vskip0.2cm
{\bf Remark 2.1} Noting that the Reilly inequality (\ref{0.1}) has been strengthened to \cite{wx1}
\be \label{r1}
\la_1\leq \fr{(n-1)|M|}{n|\om|}\left(\fr{\omega_n}{|\om|}\right)^{1/n}
\en
with equality holding if and only $M$ is a hypersphere, we believe that a stronger form of (\ref{0.2}) is valid.
\begin{conjecture} If the conditions are as in Theorem \ref{th1}, then
\be\label{c1}
\sum_{i=1}^n \la_i\leq \fr{(n-1)|M|}{|\om|}\left(\fr{\omega_n}{|\om|}\right)^{1/n}.
\en
Moreover, the equality holds in (\ref{c1})  if and only $M$ is a round sphere.
\end{conjecture}
\vskip0.2cm
{\bf Remark 2.2}  It is easy to see from (\ref{pth1.8}) that
\be \label{r2.2}
\la_n\leq\fr{|\pa\om|^2}{4|\om|^2}
\en
which is also new. It would be interesting to know the best possible upper bound for $\la_n$.

\section{Proofs of Theorems \ref{th3} and \ref{th4}}
\setcounter{equation}{0}
In this section, we shall prove  Theorems \ref{th3} and \ref{th4}. Before doing this,  let us recall some known facts. Let $\{\phi_i\}_{i=0}^{\infty}$ be  orthonormal eigenfunctions corresponding to the eigenvalues $\{\xi_i\}_{i=0}^{\infty}$ of the problem (\ref{0.4}). That is,
\be\label{.}
 \left\{\begin{array}{l}
 \Delta^2 \phi_i= 0  \ \ {\rm in \ \ } \om, \\
\pa_{\nu}\phi_i = \pa_{\nu}(\D\phi_i)+\xi_i\phi_i =0 \ \ {\rm on \ \ } \pa\om
\\
\int_{\pa\om}\phi_i\phi_j=\delta_{ij}.
\end{array}\right.
\en
For each $k=1,\cdots,$ the variational characterization for $\xi_k$  is given by
\be\label{t3.1}
\xi_k=\underset{\underset{\int_{\pa\om}\phi\phi_i=0, i=0, 1,\cdots, k-1}{\phi\in H^2(\om), \pa_{\nu}\phi|_{\pa\om}=0, \phi|_{\pa\om}\neq 0}}{\inf}\fr{\int_{\om}(\D\phi)^2}{\int_{\pa\om}\phi^2}.
\en
Let $\om$ be a bounded domain in $\R^n$ and  $\om^*$  the ball centered at the origin in $\R^n$ such that $|\om^*|=|\om|$. The moments of inertia of  $\om$ with respect to the hyperplanes $x_k=0$, are defined as
\be\label{3.1}
J_k(\om)=\int_{\om}x_k^2 \ {\rm for \ all\ } k\in\{1,\cdots,n\}.
\en
By summation over $k$, we obtain the polar moment of inertia of $\om$ with
respect to the origin denoted by
\be\label{3.2}
J_0(\om)=\sum_{k=1}^n\int_{\om}x_k^2.
\en
Note that $J_0(\om)$ depends on the position of the origin.
 In fact, $J_0(\om)$ is smallest
when the origin coincides with the center of mass of $\om$, i.e. when we have
\be
\int_{\om} x_k=0, \ k=1,\cdots,n.
\en
We need the following well known  isoperimetric property \cite{bf}, \cite{hps}:
\begin{theorem}\label{t03.1}  Among all domains $\om$ of prescribed $n$-volume, the ball $\om^*$
centered at the origin has the smallest polar moment of inertia, that is,
\be
J_0(\om)\geq J_0(\om^*), \ \ \om\in \mathcal{O},
\en
for all bounded domain $\om$ of prescribed $n$-volume $|\om|$, with equality if and only if
$\om$ coincides with $\om^*$.
\end{theorem}
By multiplication over $k$ in (\ref{3.1}), we obtain a quantity denoted by $J(\om)$,
\be
J(\om)=\overset{n}{\underset{k=1}{\prod}}J_k(\om)
 \en
 which satisfies the
 following isoperimetric inequality  \cite{bl}, \cite{he2}:
\be\label{4.2}
J(\om)\geq J(\om^*)=\fr{|\om|^{n+2}}{(n+2)^n\omega_n^2}
\en
with equality if and only if $\om$ is an ellipsoid symmetric with respect to the hyperplanes $x_k=0, k=1,\cdots, n.$
\vskip0.3cm
{\it Proof of Theorem \ref{th3}.} By a translation of the origin  in $\R^n$, we can assume that
\be\label{t3.1}
\int_{\Om} x_i=0, \ i=1,\cdots,n.
\en
For each $i\in \{1,\cdots, n\}$, let $g_i$ be the solution of the problem
\be\label{pth3.3}
 \left\{\begin{array}{l}
  \D g_i= x_i\ \ \ \ \ \mbox{in}\ \ \  \Om, \\
  \pa_{\nu}g_i|_{\pa \Om}=0,\\
  \int_{\pa \Om} g_i=0.
\end{array}\right.
\en
We {\it claim} that if the
coordinate functions $x_1,\cdots,x_n$ are chosen properly,
then
\be
g_i\bot{\rm span}\{\phi_0,\cdots, \phi_{i-1}\}, \ i=1,\cdots,n.
\en
To see this, let us fix a set of coordinate functions $x_1,\cdots, x_n$ and the solutions $g_1,\cdots,g_n$
as above.
Consider the $n\times n$ matrix $H=(h_{ji})$ with
$h_{ji}=\int_{M} g_j\phi_i$, for $i,j=1,2,\cdots,n.$ One can find  an upper triangle matrix $S=(s_{ji})$ and an
orthogonal matrix $T=(t_{ji})$ such that $S=TH$, that is,
\be\label{pth3.4}
s_{ji}=\sum_{k=1}^n t_{jk}h_{ki}=\int_M \sum_{k=1}^n  t_{jk}g_k  \phi_i
=0,\ \  1\leq i<j\leq n.
\en
Letting $y_j=\sum_{k=1}^n  t_{jk}x_k$, $\tilde{g}_j=\sum_{k=1}^n  t_{jk}g_k$,
we have from (\ref{pth3.3}) and (\ref{pth3.4}) that
\be\label{pth3.5}
\left\{\begin{array}{l}
  \D\tilde{g}_i= y_i\ \ \ \ \ \mbox{in}\ \ \  \Om, \\
  \pa_{\nu}\tilde{g}_i|_{\pa \Om}=0,\\
  \int_{\pa \Om} \tilde{g}_i=0
\end{array}\right.
\en
and
\be
\tilde{g}_i\bot{\rm span}\{\phi_0,\cdots, \phi_{i-1}\}, \ i=1,\cdots,n.
\en
Since $T=(t_{ji})$ is an orthogonal matrix, $y_1,\cdots,y_n$ are also coordinate functions of $\R^n$. Thus, our {\it claim} is true. Denoting these coordinate functions and the solutions of  (\ref{pth3.5}) still by $x_1, x_2,···, x_n,$ and $g_1,\cdots,g_n$, respectively, we conclude from (\ref{t3.1}) that
\be\label{pth3.6}
\xi_i\leq \fr{\int_{\Om} x_i^2}{\int_{\pa \Om}  g_i^2}, \ \ i=1,\cdots,n.
\en
From divergence theorem we know that
\be
\int_{\Om} x_i^2 =\int_{\Om} x_i \D g_i=-\int_{\Om}\langle\na x_i, \na g_i\ra=-\int_{\pa \Om} g_i\pa_{\nu}x_i,
\en
which gives
\be\label{pth3.7}
\left(\int_{\Om} x_i^2\right)^2\leq\int_{\pa\Om} (\pa_{\nu} x_i)^2\int_{\pa\Om} g_i^2,\ i=1,\cdots,n.
\en
Combining (\ref{pth3.6}) into  (\ref{pth3.7}), we infer
\be\label{pth3.8}
\xi_i{\int_{\Om} x_i^2}\leq\int_{\pa \Om}(\pa_{\nu} x_i)^2,\ i=1,\cdots,n,
\en
which implies that
\be\label{pth3.9}
\sum_{i=1}^n\fr 1{\xi_i}\geq\sum_{i=1}^n\fr{\int_{\Om} x_i^2}{\int_{\pa \Om}(\pa_{\nu} x_i)^2}.
\en
Using the arithmetic-geometric mean inequality and the isoperimetric inequality (\ref{4.2}), we have
\be\label{pth3.10}\no
\sum_{i=1}^n\fr{\int_{\Om} x_i^2}{\int_{\pa \Om}(\pa_{\nu} x_i)^2}&\geq&\fr{n\left(\overset{n}{\underset{j=1}{\prod}}\int_{\Om} x_i^2\right)^{1/n}}{\left(\overset{n}{\underset{j=1}{\prod}}\int_{\pa \Om}(\pa_{\nu} x_i)^2\right)^{1/n}}\\  \no
&\geq&\fr{n|\om|^{1+\fr 2n}}{(n+2)\omega_n^{2/n}}\cdot\fr 1{\left(\overset{n}{\underset{j=1}{\prod}}\int_{\pa \Om}(\pa_{\nu} x_i)^2\right)^{1/n}}\\ \no
&\geq&\fr{n|\om|^{1+\fr 2n}}{(n+2)\omega_n^{2/n}}\cdot\fr 1{\fr 1n \sum_{i=1}^n\int_{\pa \Om}(\pa_{\nu} x_i)^2}\\
&=&\fr{n^2|\om|^{1+\fr 2n}}{(n+2)|\pa\om|\omega_n^{2/n}}
\en
with equality holding if and only if
\be\label{pth3.11}& &
\fr{\int_{\Om} x_1^2}{\int_{\pa \Om}(\pa_{\nu} x_1)^2}=\cdots=\fr{\int_{\Om} x_n^2}{\int_{\pa \Om}(\pa_{\nu} x_n)^2},\\ & &
\label{pth3.12}
\overset{n}{\underset{j=1}{\prod}}\int_{\Om} x_i^2=\fr{|\om|^{n+2}}{(n+2)^n\omega_n^2}
\en
and
\be\label{pth3.13}
\int_{\pa \Om}(\pa_{\nu} x_1)^2=\cdots=\int_{\pa \Om}(\pa_{\nu} x_n)^2.
\en
Combining (\ref{pth3.9}) and (\ref{pth3.10}), one gets (\ref{0.5}). If the equality holds in (\ref{0.5}),
then (\ref{pth3.6}), (\ref{pth3.7}), (\ref{pth3.8}), (\ref{pth3.9}) and (\ref{pth3.10}) should take equality. It follows that
\be\label{pth3.14}
\int_{\Om} x_1^2=\int_{\Om} x_2^2=\cdots \int_{\Om} x_n^2=\fr{|\om|}{n+2}\cdot\left(\fr{|\om|}{\omega_n}\right)^{2/n}
\en
and so
\be\label{pth3.15}
\int_{\Om}\sum_{i=1}^n x_i^2=\fr{n|\om|}{n+2}\cdot\left(\fr{|\om|}{\omega_n}\right)^{2/n}.
\en
Consequently, we conclude from Theorem \ref{t03.1}  that $\om$ is a ball. On the other hand, if $\om$ is a ball of radius $R$ in $\R^n$, then
\be
\sum_{i=1}^n\fr 1{\xi_i}=\fr n{\xi_1}=\fr n{\fr{(n+2)}{R^3}}=\fr{n^2|\om|\left(\fr{|\om|}{\omega_n}\right)^{2/n}}{(n+2)|\pa\om|}.
\en
This completes the proof of Theorem \ref{th3}.\hfill $\Box$
\vskip0.3cm
{\bf Remark 3.1.} Consider a more general eigenvalue problem :
\be\label{r3.2}
 \left\{\begin{array}{l}
 \Delta^2 u= 0  \ \ {\rm in \ \ } \om, \\
\pa_{\nu}u = \pa_{\nu}(\D u)+\zeta\rho u =0 \ \ {\rm on \ \ } \pa\om,
\end{array}\right.
\en
where $\rho$ is a continuous positive function on $\pa\om$.  The eigenvalues of this problem  can be arranged as (counted with multiplicity):
\be
0=\zeta_0<\zeta_1\leq\zeta_2\leq\cdots\nearrow +\infty.
\en
When $\om $ is a bounded domain with smooth boundary in $\R^n$,  one can use similar arguments as in the proof of (\ref{pth3.17}) to show that the first $n$ nonzero eigenvalues of the
problem (\ref{r3.2}) satisfy
\be\label{pth3.17}
\overset{n}{\underset{j=1}{\prod}}\xi_j\leq
\left(\fr{\omega_n}{|\om|}\right)^2\cdot\left(\fr{(n+2)}{n|\om|}\int_{\pa\om}\fr 1{\rho}\right)^n,
\en
with equality holding implies that $\om$ is an ellipsoid. To see this, let us take an orthonormal set of
eigenfunctions $\{\psi_i\}_{i=0}^{\infty}$ corresponding to the
eigenvalues $\{\zeta_i\}_{i=0}^{\infty}$, that is,
\be\label{r3.4}
 \left\{\begin{array}{l}
 \Delta^2 \psi_i= 0  \ \ {\rm in \ \ } \om, \\
\pa_{\nu}\psi_i = \pa_{\nu}(\D\psi_i)+\zeta_i\rho\psi_i =0 \ \ {\rm on \ \ } \pa\om.\\
\int_{\pa\om}\rho\psi_i\psi_j=\delta_{ij}
\end{array}\right.
\en
The variational characterization for $\zeta_k$  is given by
\be\label{4.3}
\zeta_k=\underset{\underset{\int_{\pa\om}\rho\psi\psi_i=0, i=0, 1,\cdots, k-1}{\psi\in H^2(\om), \pa_{\nu}\psi|_{\pa\om}=0, \psi|_{\pa\om}\neq 0}}{\inf}\fr{\int_{\om}(\D\psi)^2}{\int_{\pa\om}\rho\psi^2}, \ k=1,\cdots.
\en
We choose the origin in $\R^n$  so that (\ref{t3.1}) holds. For each $i\in \{1,\cdots, n\}$, let $h_i$ be the solution of the problem
\be\label{r3.3}
 \left\{\begin{array}{l}
  \D h_i= x_i\ \ \ \ \ \mbox{in}\ \ \  \Om, \\
  \pa_{\nu}h_i|_{\pa \Om}=0,\\
  \int_{\pa \Om}\rho h_i=0.
\end{array}\right.
\en
As in the proof of Theorem \ref{th3}, we can assume that
\be\label{r3.4}
\int_{\pa\om}\rho x_i\psi_j=0, \ i=1, 2, \cdots,n, \  j<i.
\en
It  follows from (\ref{4.3}) that
\be\label{r3.5}
\zeta_i\int_{\pa\om}\rho h_i^2\leq \int_{\om}x_i^2, \ i=1,\cdots,n.
\en
Since
\be\label{r3.6}
\int_{\om} x_i^2=\int_{\om}x_i\Delta h_i=-\int_{\pa\om}h_i\pa_{\nu}x_i\leq\left(\int_{\pa\om}\rho h_i^2\right)^{1/2}\left(\int_{\pa\om}\fr{(\pa_{\nu}x_i)^2}{\rho}\right)^{1/2},
\en
we infer from (\ref{r3.5}) that
\be\label{r3.6}
\zeta_i\int_{\om} x_i^2\leq \int_{\pa\om}\fr{(\pa_{\nu}x_i)^2}{\rho}, \ i=1,\cdots,n.
\en
By multiplication over $i$,  one gets
\be\label{pth3.16}\no
\overset{n}{\underset{j=1}{\prod}}\zeta_j\cdot \overset{n}{\underset{i=1}{\prod}}\int_{\Om}x_i^2
&\leq& \overset{n}{\underset{i=1}{\prod}}\int_{\pa \Om}\fr{(\pa_{\nu} x_i)^2}{\rho}\\ \no
&\leq& \left(\fr 1n\sum_{i=1}^n \int_{\pa \Om}\fr{(\pa_{\nu} x_i)^2}{\rho}\right)^n
\\ &=&\left(\fr 1n\int_{\pa\om}\fr 1{\rho}\right)^n,
\en
which, combining with (\ref{4.2}), yields (\ref{pth3.17}). Also, when the equality holds (\ref{pth3.17}), we must have the equality in (\ref{4.2}) and so   $\om$ is an ellipsoid.
\vskip0.3cm
{\it Proof of Theorem \ref{th4}} Let $u_0, u_1, u_2, \cdots, $ be  orthonormal eigenfunctions corresponding to the eigenvalues
$0, \lambda_{1,\tau}, \lambda_{1,\tau}, \cdots,$  that is,
\begin{eqnarray*}
\left\{\begin{array}{ccc} \D^2 u_i-\tau\D u_i
=0,&&~\mbox{in}~~\Omega, \\[2mm]
\frac{\p^2u_i}{\p\nu^2}=0, &&~~\mbox{on}~~\partial \Omega,\\[2mm]
\tau\frac{\p u_i}{\p \nu}-\mathrm{div}_{\p \Omega}\(\nabla^2 u_i(\nu)\)-\frac{\p\D u_i}{\p\nu}=-\lambda_{i,\tau} u_i,&&~~\mbox{on}~~\partial \Omega,\\[2mm]
\inpo u_i u_j =\delta_{ij}.
\end{array}\right.
\end{eqnarray*}
Note that $u_0= 1/\sqrt{|\pa \om|}$. Using the same discussions as in the proof of Theorem \ref{th1}, we can assume   that
\be \int_{\pa\om} x_i u_j=0, \ i=1,\cdots, n, \ j=0,\cdots, i-1.\en
Thus, we have from (\ref{0.7}) that
\begin{eqnarray}\label{c4}
\lambda_{i,\tau}\int_{\pa\om} x_i^2\leq \int_{\om}
\(|\na^2x_i|^2 +\tau|\na x_i|^2\) = \tau |\om|, \ i=1,\cdots,n.
\end{eqnarray}
As in the proof of Theorem \ref{th1}, we have for each $i\in\{1,\cdots,n\}$ that
\be\label{pt4.1}
|\om|^2=\left(\int_{\pa\om} x_i\pa_{\nu}x_i\right)^2\leq \left(\int_{\pa\om} x_i^2\right)\left(\int_{\pa\om} (\pa_{\nu}x_i)^2\right)
\en
with equality holding if and only $\pa_{\nu}x_i=\eta_i x_i$ for some constant $\eta_i\neq 0$.

Multiplying (\ref{c4}) by $\int_{\pa\om} (\pa_{\nu}x_i)^2$ and using (\ref{pt4.1}), we get
\begin{eqnarray}\label{c5}
\lambda_{i,\tau}|\om|^2\leq \tau|\om|\int_{\pa\om} (\pa_{\nu}x_i)^2, \ 1,\cdots, n.
\en
 Dividing by $|\om|^2$ and summing over $i$, one gets
\begin{eqnarray}\label{pt4.2}
\sum_{j=1}^{n}\lambda_{j,\tau}\leq
\frac{\tau}{|\om|}\int_{\pa\om}\sum_{i=1}^n (\pa_{\nu}x_i)^2=\frac{\tau|\pa\om|}{|\om|}.
\en
This proves (\ref{th4.1}). Moreover, if equality holds in (\ref{th4.1}), then
$\pa_{\nu}x_i=\eta_i x_i, $ for some nonzero constants $\eta_i\  i=1,\cdots, n.$
It follows that
\be
\sum_{i=1}^n \eta_i^2 x_i^2=1 \ \ {\rm on\ } \pa\om.
\en
If $z=\sum_{i=1}^n \eta_i^2 x_i^2$, then the outward unit normal of $\pa\om$ is given by
\be \label{pt4.3}
\nu=\fr{\na z}{|\na z|}.
\en
Note that
\be \label{pt4.4}
\nu=(\pa_{\nu}x_1,\cdots, \pa_{\nu}x_n)=(\eta_1 x_1,\cdots, \eta_n x_n).
\en
Comparing (\ref{pt4.3}) and (\ref{pt4.4}), we infer $\eta_1=\eta_2=\cdots=\eta_n$, which shows that
$\pa\om$ is a hypersphere and so $\om$ is a ball. On the other hand, we have
\be
\lambda_{1,\tau}(\mathbf{B}_R^n)=\cdots\lambda_{n,\tau}(\mathbf{B}_R^n)=\fr{\tau}{R}=\fr{\tau|\pa\mathbf{B}_R^n|}
{n|\mathbf{B}_R^n|},
\en
that is, the equality holds for balls in (\ref{th4.1}). \hfill $\Box$
\vskip0.3cm
\begin{conjecture}
 Under the same assumptions of Theorem \ref{th4}, we have
\be\label{c3.1} \overset{n}{\underset{j=1}{\prod}}\lambda_{j,\tau}\leq\fr{\tau^n\omega_n}{|\om|}
\en
with equality holding if and only if $\om$ is a ball.
\end{conjecture}
\vskip0.2cm
It should be mentioned that if $\om$ is convex, then the above conjecture is true. To see this, it is enough to  take the products of the $n$ inequalities in (\ref{c4}) and use Lemma \ref{l4.1}.

\section{Proofs of Theorems \ref{th7} and \ref{th8}}
\setcounter{equation}{0}
In this section, we prove Theorems \ref{th7} and \ref{th8}. We shall need the following result
\cite{hps}.
\vskip0.3cm
\begin{lemma} \label{l4.1}
 Let $\om$ be a bounded convex domain in $\R^n$. Assume that   origin coincides with the center of mass of $\pa\om$, that is,
\be\label{pth7.1} \int_{\pa\om}x_i ds=0,\ i=1,\cdots,n.
\en
Then we have
\be \label{pth7.2}
\overset{n}{\underset{i=1}{\prod}}\int_{\pa\om} x_i^2 ds\geq \overset{n}{\underset{i=1}{\prod}}\int_{\pa\om^*} x_i^2 ds=\fr{|\om|^{n+1}}{\omega_n},
\en
with equality if and only if $\om=\om^*$.
\end{lemma}
\vskip0.3cm
We prove the following result from which  Theorem \ref{th7} follows.

\begin{theorem} \label{t4.1}  Let $\beta\geq 0$ and $\Omega$ be a bounded domain with smooth boundary $\p \Omega$ in
 $\R^n$.  Let $\rho$ be a positive continuous function on $\pa\om$ and denote by $0<\eta_{1, \beta}\leq\eta_{2, \beta}\leq\cdots\leq \eta_{n, \beta}\leq\cdots $
 the
 eigenvalues of the  problem :
\be\label{t7.1}\left\{\begin{array}{l} \Delta u =0 \ \ \ \ \ \
\ \ \ \ \ \ \ \ \ \ \ \ \ \ \ \ \ {\rm in  \ }\ \om,\\ -\beta\ov{\Delta}
u+\pa_{\nu} u= \eta\rho u\ \ \ \ \quad ~~ {\rm on  \ } \pa \om,
\end{array}\right.
\en
Then we have
\be\label{t7.2}
\sum_{i=1}^n \eta_{i,\b}\leq \fr 1{|\om|^2}\left((|\om|+\beta|\pa\om|)\int_{\pa\om}\rho^{-1}
-\fr{\beta}n\left(\int_{\pa\om}\fr 1{\sqrt{\rho}}\right)^2\right).
 \en
Furthermore, if $\rho$ is constant, the equality holds in (\ref{t7.2}) if and only if $\om$
is a ball.
\end{theorem}

{\it Proof of Theorem \ref{t4.1}.} Let $u_0, u_1, u_2,\cdots$ be
 orthonormal  eigenfunctions
corresponding to the eigenvalues
$0=\eta_{0,\b}<\eta_{1,\b}\leq \eta_{2,\b}\leq \cdots, $ of
the problem (\ref{t7.1}), that is,
\be\label{t7.3}\left\{\begin{array}{l} \Delta u_i =0 \ \ \ \
\ \ \ \ \ \ \ \ \ \ \ \ \ \ \ \ \ \ \ \ {\rm in  \ }\ \om,\\
-\beta\ov{\Delta}
u+\pa_{\nu} u_i= \eta_i\rho u_i\ \ \ \ \quad  {\rm on  \ } \pa \om,\\
\int_{\pa \om} \rho u_iu_j=\delta_{ij}.
\end{array}\right.
\en
Note that $u_0$ is a constant function $1/(\int_{\pa\om}\rho)^{1/2}$. The eigenvalues $\eta_{i,\b}, i=1, 2,\cdots,$ are
characterized by
\be
\label{t7.4}
\eta_{i,\b}=\underset{\underset{\int_{\pa\om}\rho u u_j=0, \ j=0,1,\cdots, i-1}{u\in H(\om)\setminus\{0\}}}{\min}
\fr{\int_{\om}|\na u|^2 +\beta\int_{\pa\om} |\ov{\na}
u|^2}{\int_{\pa\om}\rho u^2}.
\en
As in the proof of Theorem \ref{th1}, we can choose the coordinate functions $x_1,\cdots,x_n$ of
$\R^n$ so that
\be
\int_{\pa\om}\rho x_iu_j=0,\ j<i, i=1,\cdots,n.
\en
Hence
\be
\eta_{i,\b}\int_{\pa\om} \rho x_i^2 &\leq& \int_{\om}|\na x_i|^2 +\beta\int_{\pa\om} |\ov{\na}
x_i|^2\\ \no &=&|\om|+\beta\int_{\pa\om} |\ov{\na}
x_i|^2\\ \label{pth7.6} &=& |\om|+\beta\left(|\pa\om|-\int_{\pa\om}(\pa_{\nu}x_i)^2\right), \ i=1,\cdots,n.
\en
We have from (\ref{pt4.1}) that
\be\label{t7.5}
|\om|^2\leq \left(\int_{\pa\om}\rho x_i^2\right)\left(\int_{\pa\om}\rho^{-1}(\pa_{\nu}x_i)^2\right).
\en
Multiplying (\ref{pth7.6}) by $\int_{\pa\om}\rho^{-1}(\pa_{\nu}x_i)^2$ and using (\ref{t7.5}), we get
\be\no
\eta_{i,\b}|\om|^2&\leq& (|\om|+\beta|\pa\om|)\int_{\pa\om}\rho^{-1}(\pa_{\nu}x_i)^2
-\beta\left(\int_{\pa\om}(\pa_{\nu}x_i)^2\right)\left(\int_{\pa\om}\rho^{-1}(\pa_{\nu}x_i)^2\right)
\\ \label{pth7.7} &\leq& (|\om|+\beta|\pa\om|)\int_{\pa\om}\rho^{-1}(\pa_{\nu}x_i)^2
-\beta\left(\int_{\pa\om}\fr 1{\sqrt{\rho}}(\pa_{\nu}x_i)^2\right)^2.
\en
Summing over $i$ and using Cauchy-Schwarz inequality, one has
\be\no
|\om|^2\sum_{i=1}^n\eta_{i,\b}
&\leq& (|\om|+\beta|\pa\om|)\int_{\pa\om}\rho^{-1}-\fr{\beta}n\left(\sum_{i=1}^n\int_{\pa\om}\fr 1{\sqrt{\rho}}(\pa_{\nu}x_i)^2\right)^2
\\ \label{pth7.8} &=&(|\om|+\beta|\pa\om|)\int_{\pa\om}\rho^{-1}
-\fr{\beta}n\left(\int_{\pa\om}\fr 1{\sqrt{\rho}}\right)^2.
\en
Dividing by $|\om|^2$, we get (\ref{th7.2}). Moreover, when $\rho$ is  constant,  the equality holds in (\ref{th7.2}) if and only if $\om$ is a ball. \hfill $\Box$

\vskip0.3cm
{\it Proof of Theorem \ref{th8}.} Let us choose the origin in $\R^n$  as the center of mass of $\pa\om$. Taking $\rho=1$ and using the same arguments as  in the proof of Theorem \ref{th7}, we can get
\be \lambda_{i,\b}\int_{\pa\om}  x_i^2 \leq |\om|+\beta\left(|\pa\om|-\int_{\pa\om}(\pa_{\nu}x_i)^2\right), \ i=1,\cdots,n.
\en
By multiplication of these inequalities, one infers
\be\no
& & \overset{n}{\underset{i=1}{\prod}}\lambda_{i,\b}\overset{n}{\underset{j=1}{\prod}}\int_{\pa\om} x_j^2\\ \no &\leq&
\overset{n}{\underset{i=1}{\prod}}\left((|\om|+\beta|\pa\om|)-\beta\int_{\pa\om}(\pa_{\nu}x_i)^2\right)
\\ \no &\leq&\left(\fr 1n\sum_{i=1}^n\left((|\om|+\beta|\pa\om|)-\beta\int_{\pa\om}(\pa_{\nu}x_i)^2\right)\right)^n
\\ \label{pth8.1} &=&\left(|\om|+\fr{(n-1)\beta}n|\pa\om|\right)^n
\en
Substituting (\ref{pth7.2}) into (\ref{pth8.1}), we obtain (\ref{th8.1}). It is clear from the proof that
equality holds in (\ref{th8.1}) if and only $\om$ is a ball. \hfill
$\Box$

\vskip0.2cm
{\bf Remark 4.1} We believe that the convexity assumption in Theorem \ref{th8} is unnecessary.

\vskip0.4cm

Department of Mathematics,
Southern University of Science and Technology

Shenzhen, 518055, GUANDONG,
P. R.  CHINA

fuquan\_fang@yahoo.com

\vskip0.6cm

xiachangyu666@163.com

\end{document}